%
%
%
%
%
\font\tenbf=cmbx10

\font\eightrm=cmr8
\font\eightit=cmti8
\font\germ=eufm10

\def\s{\hbox{\germ S}}
\def\sectiontitle#1\par{\vskip0pt plus.1\vsize\penalty-250
\vskip0pt plus-.1\vsize\bigskip\vskip\parskip
\message{#1}\leftline{\tenbf#1}\nobreak\vglue 5pt}
\def\wt{\widetilde}
\def\wh{\widehat}
\def\ds{\displaystyle}
\def\eno{\eqalignno}
\def\ld{\lambda}

\def\al{\alpha}

\def\frac#1#2{{#1\over#2}}

\magnification=\magstep1
\vsize=8.6truein 
\parindent=15pt
\nopagenumbers
\baselineskip=10pt
\line{\hfil\it
Dedicated to the memory of }
\baselineskip=12pt
\line{\hfil\it 
Marcel--Paul Sch\"utzenberger
}
\vglue 2pc 
\baselineskip=13pt
\headline{\ifnum\pageno=1\hfil\else%
{\ifodd\pageno\rightheadline \else \leftheadline\fi}\fi}
\def\rightheadline{\hfil\eightit 
Quantum Schubert polynomials and quantum Schur functions
\quad\eightrm\folio}
\def\leftheadline{\eightrm\folio\quad 
\eightit 
Anatol N. Kirillov 
\hfil}
\voffset=2\baselineskip 
\centerline{\tenbf QUANTUM \hskip 0.1cm SCHUBERT \hskip 0.1cm
POLYNOMIALS\hskip 0.1cm AND 
}
\vglue 3pt
\centerline{\tenbf QUANTUM
\hskip 0.1cm SCHUR  \hskip 0.1cm FUNCTIONS 
}
\vglue 10pt
\centerline{\eightrm 
ANATOL N. KIRILLOV
}
\baselineskip=18pt
\centerline{\eightit 
CRM, University of Montreal
}
\baselineskip=10pt
\centerline{\eightit 
C.P. 6128, Succursale A, Montreal (Quebec) H3C 3J7, Canada
}
\baselineskip=12pt
\centerline{\eightit
and }
\baselineskip=12pt
\centerline{\eightit 
Steklov Mathematical Institute,
}
\baselineskip=10pt
\centerline{\eightit 
Fontanka 27, St.Petersburg, 191011, Russia
}
\vglue 20pt

\centerline{\eightrm ABSTRACT}
{\rightskip=1.5pc
\leftskip=1.5pc
\eightrm\parindent=1pc 
We introduce the quantum multi--Schur functions, quantum factorial 
Schur functions and quantum Macdonald polynomials. We prove that for 
restricted vexillary permutations the quantum double Schubert polynomial 
coincides with some quantum multi--Schur function and prove a quantum 
analog of the N\"agelsbach--Kostka and Jacobi--Trudi formulae for 
the quantum double Schubert polynomials in the case of Grassmannian 
permutations. We prove also an analog of the 
Billey--Jockusch--Stanley formula for quantum Schubert polynomials. 
Finally we formulate two conjectures about the structure of quantum 
double and quantum Schubert polynomials for 321--avoiding permutations.
} 
\vglue12pt
\baselineskip=13pt
\overfullrule=1pt
\def\qed{\hfill$\vrule height 2.5mm width 2.5mm depth 0mm$}

{\bf \S 1. Introduction.}
\vskip 0.3cm

The cohomology ring of the flag variety $Fl_n=SL_n/B$ is isomorphic to 
the quotient ring of the polynomial ring by the ideal generated by 
symmetric polynomials without constant term. The Schubert cycles give a 
linear basis of the cohomology ring and they are represented by Schubert 
polynomials. A.~Lascoux and M.-P.~Sch\"utzenberger [LS1], [LS2] introduced the 
Schubert polynomials ${\s}_w$ as a stable, homogeneous basis of $P_n:={\bf 
Z}[x_1,\ldots ,x_n]$ indexed by permutations $w\in S^{(n)}$. We refer the 
reader to [M1] for detailed account on Schubert polynomials. It is 
well--known ([M1], (4.8)) that if $w$ is a Grassmannian permutation of 
shape $\ld$, then ${\s}_w$ is the Schur function $s_{\ld}(X_r)$, where $r$ 
is the unique descent of $w$ and $X_r=x_1+\cdots +x_r$. More generally, 
if $w$ is a vexillary permutation ([LS1]; [M1], Chapter I) with shape $\ld 
=(\ld_1,\ldots ,\ld_m)$ and flag $\phi =(\phi_1,\ldots ,\phi_m)$ ,
then ${\s}_w$ is a multi--Schur function ([LS1], Section~3; [M1], 
Chapter~II), namely
$${\s}_w=s_{\ld}(X_{\phi_1},\ldots ,X_{\phi_m}).\eqno (1.1)
$$

The main goal of this note is 

$\bullet$ To introduce the quantum multi--Schur functions  and prove an 
analog of the formula (1.1) for the double quantum Schubert polynomials 
which corresponds to the restricted vexillary permutations 
(see Definition~7).

$\bullet$ To prove the N\"agelsbach--Kostka and Jacobi--Trudi type 
formulae for the quantization of Schur function 
$s_{\ld}(X_r)$, $\ld\subset ((n-r)^r)$, and more generally, for the 
quantum double Schubert polynomials corresponding to the Grassmannian 
permutations. 

We define a quantum Schur 
function $\wt s_{\ld}(X_r)$ as the quantum Schubert polynomial $\wt{\s}_w$ 
corresponding to the Grassmannian permutation $w\in S_n$ of shape $\ld$ 
and descent $r$. Follow to A.~Lascoux~[L], we define the quantum 
factorial Schur function $\wt s_{\ld}(X_r\| a)$ as the quantum double 
Schubert polynomial $\wt{\s}_w(x,a)$ corresponding to the same 
Grassmannian permutation $w$. We introduce also quantum column--flagged 
Schur functions and study the quantum double Schubert polynomials for 
321--avoiding permutations.

For the reader's convenience, we formulate below our main results:
\vskip 0.2cm

$\bullet$  Let $\ld\subset ((n-r)^r)$ be a partition, then
$$\wt s_{\ld}(X_r)=\det (\wt e_{\ld'_i-i+j}(X_{r-1+j}))_{1\le i,j\le n-r},
$$
where $\wt e_k(x):=e_k(x|q)$ is the quantum elementary symmetric 
polynomial of degree $k$ (Theorem~1).

$\bullet$ Let $w\in S_n$ be a dominant permutation of shape $\ld$. 
Then
$$\wt{\s}_w:=\wt x^{\ld}=\det\left(\wt h_{\ld_i-i+j}(X_i)\right)_{1\le 
i,j\le l(\ld )},
$$
where $\wt h_k(X_r)=\det\left(\wt e_{1-i+j}(X_{r-1+j})\right)_{1\le 
i,j\le k}$ is the quantum complete homogeneous symmetric polynomial of 
degree $k$ in the variables $X_r=x_1+\cdots +x_r$, and $\wt x^{\ld}$ is 
the quantization (see Section (2.3) below) of monomial $x^{\ld}$ 
(Corollary~3; for more general results, see Corollaries~4 and 6).

$\bullet$  Let $w$ be a restricted vexillary permutation with shape 
$\ld =\ld (w)$ and flag $\theta =\theta (w)$ and 
let $\ld'=\ld (w^{-1})$ and $\wt\theta =\theta (w^{-1})$ be shape and 
flag of the inverse permutation $w^{-1}$. Then
$$\wt{\s}_w(x,y)=s_{\ld}^q(Z_1,Z_2,\ldots ,Z_m), 
$$
where $Z_i=X_{\theta_i}-Y_{\wt\theta_{\ld'_i}}$, $1\le i\le m=l(\ld )$, and
$s^q_{\ld}(Z_1,\ldots ,Z_m)$ is the quantum row--flagged Schur function 
(Theorem~5).

$\bullet$  Let $w\in S_n$ be a Grassmannian permutation with 
shape $\ld$ and descent $r$, and let $\wt\theta =\theta (w^{-1})$ 
be the flag of inverse permutation $w^{-1}$. Then
$$\wt{\s}_w(x,y)=\det\left(\wt e_{\ld'_i-i+j}(X_{r-1+j}-
Y_{\wt\theta_i})\right)_{1\le i,j\le n-r}, 
$$
where $\wt e_m(X_k-Y_l):=\ds\sum_{j=0}^m\wt e_{m-j}(X_k)h_j(Y_l)$ ~(Theorem~6).
 
In Section~(3.4) we introduce the quantum Macdonald polynomials $\wt 
P_{\ld}(X_r;p,t)$ and prove the quantum Cauchy identity for them.
\vskip 0.3cm

{\bf Acknowledgments.} The author would like to thank N.~Bergeron and 
N.A.~Liskova for fruitful discussions. This work was initiated during my 
stay at the University of Tokyo (1995/1996) and was completed at the CRM,
University of Montreal. I would like 
to thank all my colleagues from the Tokyo and Montreal Universities for  
very fruitful working atmosphere.

\vskip 0.5cm
{\bf \S 2. Quantization map.}
\vskip 0.3cm

In this section we review some basic properties of the 
quantization map. Originally, construction of quantization map appeared in 
[FGP]; independently, construction of quantization map was introduced in 
[KM] in a different form, using the Interpolation formula (see [M1], (6.8)) 
and quantum double Schubert polynomials (see [KM], Section~5). It can be 
shown that two forms of quantization mentioned above are equivalent. Let 
us remind a construction of quantization map from [KM].

\vskip 0.3cm
{\bf (2.1)} Quantum double Schubert polynomials.
\vskip 0.2cm

Let $x=(x_1,\ldots ,x_n)$, $y=(y_1,\ldots ,y_n)$ be two sets of 
variables, put
$$\wt{\s}_{w_0}(x,y):=\wt{\s}_{w_0}^{(q)}(x,y)=\prod_{i=1}^{n-1}\Delta_i
(y_{n-i}~|~x_1,\ldots ,x_i), \eqno (2.1)
$$
where $w_0$ is the longest element of the symmetric group $S_n$, and 
~~$\Delta_k(t~|~x_1,\ldots ,x_k):=\ds\sum_{i=0}^ke_i(x|q)t^{k-i}=$
$$ =\det\pmatrix{ x_1+t & q_1&0 &\ldots &\ldots 
&\ldots &0\cr
-1 & x_2+t & q_2 &0 &\ldots &\ldots & 0\cr
0 & -1 &x_3+t & q_3 & 0 &\ldots & 0 \cr
\vdots &\ddots &\ddots & \ddots &\ddots &\ddots &\vdots \cr
0&\ldots & 0 &-1&x_{k-2}+t &q_{k-2} & 0 \cr
0 &\ldots &\ldots & 0 &-1 & x_{k-1}+t & q_{k-1}\cr
0 &  \ldots &\ldots &\ldots & 0 & -1 & x_k+t} \eqno (2.2)
$$
\vskip 0.3cm

Polynomials $\wt e_i(x):=e_i(x|q)$ defined in the LHS(2.2)  are called 
the quantum elementary symmetric polynomials (cf. [GK], [C]). 

\vskip 0.3cm
{\bf Definition 1} ([KM], Section 3). {\it For each permutation $w\in S_n$, 
the quantum 
double Schubert polynomial $\wt{\s}_w(x,y)$ is defined to be
$$\wt{\s}_w(x,y)=\partial_{ww_0}^{(y)}\wt{\s}_{w_0}(x,y), \eqno (2.3)
$$
where divided difference operator $\partial_{ww_0}^{(y)}$ acts on the $y$ 
variables.}

\vskip 0.3cm
{\bf Definition 2} ([KM], Section~4). {\it The quantum Schubert polynomial 
$\wt{\s}_w(x)$ is defined to be 
the specialization $y_1=0,\ldots ,y_n=0$ of quantum double Schubert polynomial
${\s}_w(x,y)$:}
$$\wt{\s}_w:=\wt{\s}_w(x)=\partial_{ww_0}^{(y)}\wt{\s}_{w_0}(x,y)|_{y=0}=
\wt{\s}_w(x,0). \eqno (2.4)
$$

\vskip 0.3cm
{\bf (2.2)} Quantum Cauchy identity.
\vskip 0.2cm

The quantum Cauchy identity is the quantum analog of the Cauchy formula in 
the theory of Schubert polynomials ([M1], (5.10)). As it shown in [KM], 
Section~4, the quantum Cauchy identity is equivalent to the statement that 
quantum Schubert polynomials $\wt{\s}_w$ are orthogonal with respect to the
quantum residue pairing.

\vskip 0.3cm
{\bf Proposition 1} ([KM], (12)).
$$\sum_{w\in S_n}\wt{\s}_w(x){\s}_{ww_0}(y)=\wt{\s}_{w_0}(x,y),\eqno (2.5).
$$
\vskip 0.2cm

{\bf Corollary 1} ([KM], (13)-(14)).
$$\sum_{w\in S_n}\wt{\s}_w(x,z){\s}_{ww_0}(y,-z)=\wt{\s}_{w_0}(x,y),\eqno (2.6)
$$
$$
\sum_{\matrix{{\scriptstyle u\in S_n,} \cr
{\scriptstyle l(u)+l(uw^{-1})=l(w)}}}\wt{\s}_u(x,z)
{\s}_{uw^{-1}}(y,-z)=\wt{\s}_w(x,y).\eqno (2.7) 
$$

\vskip 0.3cm
{\bf (2.3)} Quantization map.
\vskip 0.2cm

Let $f\in P_n={\bf Z}[x_1,\ldots ,x_n]$ be a polynomial. According to the 
Interpolation formula ([M1], (6.8)),
$$f(x)=\sum_{w\in S^{(n)}}\partial_w^{(y)}f(y){\s}_w(x,y),
$$
where $S^{(n)}$ is the set of all permutations $w\in S_{\infty}$ such 
that the code $w$ (see [M1], p.9) has length $\le n$.
We define a quantization $\wt f$ of the polynomial $f$ by the rule
$$\wt f(x)=\sum_{w\in 
S^{(n)}}\partial_w^{(y)}f(y)\wt{\s}_w(x,y)|_{_{\overline P_n}},
$$
where for a polynomial $f\in\overline P_{\infty}$, the symbol 
$f|_{\overline P_m}$ means the restriction of $f$ to the ring of 
polynomials $\overline P_m$, i.e. the specialization 
$x_{m+1}=x_{m+2}=\cdots =0$, and $q_m=q_{m+1}=\cdots =0$.

Hence, the quantization is a ${\bf Z}[q_1,\ldots ,q_{n-1}]$ -- linear map 
$P_n\to\overline P_n$.

The main property of quantization is that it preserves the pairings (see 
[KM], Section~5, for further details)
$$\langle\wt f,\wt g\rangle_Q=\langle f,g\rangle ,\ \ f,g\in P_n.
$$

\vfil\eject
{\bf \S 3. Quantum Schur functions.}

\vskip 0.3cm
{\bf Definition 3.} {\it Let $\ld$ be a partition, 
$\ld\subset ((n-r)^r)$, $1\le r<n$. The quantum Schur function 
$\wt s_{\ld}(X_r)$ is 
defined to be the quantum Schubert polynomial $\wt{\s}_w$, corresponding
to the Grassmannian permutation $w\in S_n$ of shape $\ld$ and descent $r$.}
\vskip 0.2cm

In other words, the quantum Schur function $\wt s_{\ld}(X_r)$ is the 
quantization of the Schur function $s_{\ld}(X_r)$. In contrast to the 
classical case $q_1=q_2=\cdots =0$, the quantum Schur functions $\wt 
s_{\ld}(X_r)$ may not be symmetric with respect to variables $X_r$.

\vskip 0.3cm
{\bf (3.1)} Quantum analog of the N\"agelsbach--Kostka formula.
\vskip 0.2cm

In this section we are going to describe the quantization $\wt s_{\ld}(X_r)$
of Schur functions $s_{\ld}(X_r)$. Our approach is based on the quantum 
Cauchy identity. Let us 
remind a few definitions from [KM].

Let $w_0$ be the longest element of the symmetric group 
$S_n$ and 
$$\wt{\s}_{w_0}(x,y):=\Delta_1(y_{n-1}~|~x_1)\Delta_2(y_{n-2}~|~x_1,x_2)
\cdots\Delta_{n-1}(y_1~|~x_1,\ldots ,x_{n-1})
$$
be the quantum double Schubert polynomial corresponding to the element 
$w_0\in S_n$. It is clear that
$$\wt{\s}_{w_0}(x,y)=\sum_{I\subset\delta_n}\wt e_I(x)y^{\delta_n-I}
=\sum_{I\subset\delta_n}\wt x^Ie_{\delta_n -I}(y),
$$
where for any $I=(i_1,\ldots ,i_n)\subset\delta_n=(n-1,\ldots ,1,0)$,
$$\wt e_I(x)=\prod_{k=1}^n\wt e_{i_k}(X_{n-k}) \eqno (3.1)
$$
is the quantum elementary polynomial (see [KM], 
Section~5.2). Here we used notation $X_m=(x_1,\ldots ,x_m)$.
\vskip 0.3cm

{\bf Theorem 1.}  {\it Let $\lambda$ be a partition, 
$\ld\subset ((n-r)^r)$ for some $1\le r<n$, and ~$\wt e_s(X_r):=e_s(X_r|q)$ 
be the quantization 
of the elementary symmetric polynomial ~$e_s(x_1,\ldots ,x_r)$ 
(see [GK], or Section (2.1)). 
Then}
$$\wt s_{\ld}(X_r)=\det (\wt e_{\ld'_i-i+j}(X_{r-1+j}))_{1\le i,j\le 
n-r}.\eqno (3.2)
$$
\vskip 0.2cm

{\bf Corollary 2.} {\it Let $\wt h_k(X_r)$ be the quantization of the 
complete homogeneous symmetric polynomial $h_k(x_1,\ldots ,x_r)$. Then}
$$\wt h_k(X_r)=\det (\wt e_{1-i+j}(X_{r-1+j}))_{1\le i,j\le k}=
\wt s_{(k)}(X_r).\eqno (3.3)
$$

{\bf Corollary 3.} {\it Let $w\in S_n$ be a dominant permutation of 
shape $\ld$. Then
$$\wt{\s}_w(x)=\wt x^{\ld}=\det\left(\wt h_{\ld_i-i+j}(X_i)\right)_{1\le 
i,j\le l(\ld )}.\eqno (3.4)
$$
Here $\wt x^{\ld}$ is the quantization of monomial 
$x^{\ld}:=x_1^{\ld_1}\cdots x_{n}^{\ld_{n}}$.}
\vskip 0.3cm

{\it Proof of the Theorem 1} (cf. [M1], Remark on p.78). Let $n=r+s$. 
We consider $S_s\times S_r$ as 
a subgroup of $S_n$, with $S_s$ permuting $1,2,\ldots ,s$ and $S_r$ 
permuting $s+1,s+2,\ldots ,s+r=n$. Let $w_0^{(s)}$ and $w_0^{(r)}$ be the 
longest elements of $S_s$ and $S_r$ respectively, and let
$u=w_0^{(s)}\times w_0^{(r)}\in S_n$ be the product of two permutations 
$w_0^{(s)}$ and $w_0^{(r)}$ ([M1], p.45, or Section~(4.5) below). Then 
$\partial_u=\partial_{w_0^{(s)}}\partial_{1^s\times w_0^{(r)}}$ and
($I\subset\delta_n$)
$$\eno{
&\partial^{(y)}_{1^s\times w_0^{(r)}}y^{\delta_n-I}=\cases{0, & if \  
$(i_{s+1},\ldots ,i_{n-1})\ne\delta_r$,\cr 1, & if \ $(i_{s+1},\ldots 
,i_{n-1})=\delta_r$,}\cr
& \cr
&\partial_{w_0^{(s)}}^{(y)}y^J=s_{J-\delta_s}(y), \ {\rm if} \ y=(y_1,\ldots 
,y_s).}
$$
Hence, 
$$\partial_u^{(y)}\wt{\s}_{w_0}(x,y)=\partial_{w_0^{(s)}}^{(y)}
\sum_{I\subset\delta_{n,r}}\wt e_I(x)y^{\delta_{n,r}-I}
=\sum_{\ld\subset ((n-r)^r)}\wt\varphi_{\ld}(x)s_{\wh\ld'}(y),\eqno (3.5)
$$ 
where $\delta_{n,r}=(n-1,n-2,\ldots ,n-r)$, and 
$$\wt\varphi_{\ld}(x):=\det (\wt e_{\ld'_i-i+j}(X_{r-1+j}))_{1\le i,j\le 
n-r}.
$$
In the formula (3.5), we used notation $\wh\ld =(\wh\ld_1,\ldots 
,\wh\ld_r)$ for the complementary partition defined by 
$\wh\ld_i=n-r-\ld_{r+1-i}$, $1\le i\le r$, and $\wh\ld'$ is the conjugate 
of $\wh\ld$.

On the other hand, repeating the arguments from [M1], (5.10), we have
$$\eno{
\partial_u^{(y)}\wt{\s}_{w_0}(x,y)=\sum_{w\in 
S_n}\wt{\s}_w(x)\partial_u^{(y)}{\s}_{ww_0}(y)&=\sum_{w\in 
Gr_{r,s}}\wt{\s}_w(x){\s}_{ww_0}(y)\cr
&=\sum_{\ld\subset ((n-r)^r)}\wt 
s_{\ld}(x)s_{\wh\ld'}(y),}
$$
where $Gr_{r,s}$ is the set of all Grassmannian permutations $w\in 
S_{r+s}$ with descent $r$.

Hence, $\wt s_{\ld}(x)=\wt\varphi_{\ld}(x)$. It remains to note, that 
$\wt s_{\ld}(x)$ is the quantization of the Schur function $s_{\ld}(X_r)$.
\qed
\vskip 0.3cm

\vskip 0.3cm
{\bf Example.} Let us take $w=13524$. Then
$$\eno{
&c(w)=(0,1,2,0), \ \ \ld (w)=(2,1);\cr
&\s_w(x)=s_{\ld}(X_3)=\left|\matrix{h_2(X_3)&h_3(X_3)\cr 1& h_1(X_3)}\right|\cr
& \cr
&=x_1^2x_2+x_1x_2^2+x_1^2x_3+x_1x_3^2+x_2^2x_3+x_2x_3^2+2x_1x_2x_3.}
$$
It follows from Theorem 1, that for quantum Schubert polynomial 
$\wt{\s}_w$ we have
$$\eno{
\wt{\s}_w&=\wt s_{\ld}(X_3)=\left|\matrix{\wt e_2(X_3)& \wt e_3(X_4)\cr
1& \wt e_1(X_4)}\right|\cr & \cr
&=s_{\ld}(X_3)+q_1(x_1+x_2)+q_2(x_2+x_3)-q_3(x_1+x_2).}
$$
On the other hand,
$$\left|\matrix{\wt h_2(X_3)& \wt h_3(X_3)\cr 1 & \wt h_1(X_3)}\right|=
s_{\ld}(X_3)+q_1(x_1+x_2)+q_2(x_2+x_3)+q_3(x_3+x_4),
$$
which is equal to $\wt s_{\ld}(X_3)$ if and only if \ $x_1+x_2+x_3+x_4=0$.
\vskip 0.3cm

{\bf (3.2)} Quantum analog of the Billey--Jockusch--Stanley formula.
\vskip 0.2cm

Let $\al =(\al_1,\ldots ,\al_n)$ be a composition such that $0\le\al_i\le 
n-i$, $i=1,\ldots ,n$. We denote by $\wt x^{\al}$ the quantization (see 
Section~(2.3)) of the monomial $x^{\al}=x_1^{\al_1}\cdots x_n^{\al_n}$.
\vskip 0.3cm

{\bf Corollary 4} (of Theorem~1). {\it Let $\al$ be a composition, 
$\al\subset\delta_n$. Then}
$$\wt x^{\al}=\det\left(\wt h_{\al_i-i+j}(X_i)\right)_{1\le i,j\le 
n}.\eqno (3.6)
$$

For definition of polynomials $\wt h_k(X_r)$, see (3.3).
\vskip 0.3cm

{\bf Corollary 5} (Billey--Jockusch--Stanley's type formula ([BJS], 
Theorem~1.1) for quantum Schubert polynomials).
{\it Let $w\in S_n$. Then
$$\wt{\s}_w=\sum_{{\bf a}\in R(w),~{\bf b}\in C({\bf a})}
\wt x^{\bf b}, \eqno (3.7)
$$
where $R(w)$ is the set of all reduced words for the permutation $w$; if 
${\bf a}\in R(w)$, then 
$C({\bf a})=\{{\bf b}=(b_1\le b_2\le\cdots\le b_p)|1\le b_i\le a_i,~
a_j<a_{j+1}\Rightarrow b_j<b_{j+1}\}$~~
is the set of all ~${\bf a}$--compatible sequences ($p=l(w)$); 
~$\wt x^{\bf b}$ is the quantization of monomial 
~~$x^{\bf b}=x_{b_1}x_{b_2}\cdots x_{b_p}$ (see Corollary~4).}

\vskip 0.3cm
{\bf (3.3)} Quantum factorial Schur functions.

\vskip 0.3cm
{\bf Definition 4.} {\it Let $\ld$ be a partition such that $\ld_1\le 
n-r$ and $l(\ld )\le r$ for some $r$. We define a quantum factorial 
Schur function  
$\wt s_{\ld}(X_r\| a)$ to be equal to the quantum double Schubert 
polynomial ~$\wt{\s}_w(x,a)$, where $w\in S_n$ is the 
Grassmannian permutation of shape $\ld$ and descent $r$. }

It follows from (2.6) that the quantum factorial Schur functions satisfy 
the following Cauchy type formula:  
$$\sum_{\ld}\wt s_{\ld}(x\| a)s_{\wh\ld'}(y\| -a)=\wt{\s}_{uw_0}(x,y),
\eqno (3.8)
$$
where $\wt{\s}_{uw_0}(x,y):=\partial_u^{(y)}\wt{\s}_{w_0}(x,y)=
\partial_{w_0^{(n-r)}}^{(y)}\left(\ds\prod_{k=1}^{n-r}
\Delta_{n-k}(y_k~|~X_{n-k})\right)$.

In Section~4.5, Theorem~6, the N\"abelsbach--Kostka type formula for the 
quantum factorial Schur functions is given.

\vskip 0.3cm
{\bf Remark.} It is easy to see that
$$uw_0=\pmatrix{1&2&\cdots &r&r+1&r+2&\cdots &n\cr 
s+1&s+2&\cdots &n&1&2&\cdots &s}\in Gr_{r,s}.
$$
Hence, it follows from Theorem~6 below that
$$\wt{\s}_{uw_0}(x,y)=\det\left(\wt e_{r-i+j}(X_{r-1+j}-Y_i)\right)_{1\le 
i,j\le n-r}.
$$

\vskip 0.3cm
{\bf Example.} Let us take $w=3412\in Gr_{2,2}$. Then $\ld (w)=(2,2)$, 
$\theta (w)=(1,2)$ and we have the following formulae for quantum 
factorial Schur function 
$$\eno{
\wt s_{22}(X_2,A_2)=\wt{\s}_w(X_2,A_2)&=\left|\matrix{\wt h_2(X_1-A_2)
&\wt h_3(X_1-A_2)\cr \wt h_1(X_2-A_2)&\wt h_2(X_2-A_2)}
\right| ~~~~~~~~ {\rm (see~Theorem~5)}\cr
& \cr
&=\left|\matrix{\wt e_2(X_2-A_1)&\wt e_3(X_3-A_1)\cr
\wt e_1(X_2-A_2)&\wt e_2(X_3-A_2)}\right| ~~~~~~~~~
{\rm (see~Theorem~6)}\cr
& \cr
&=q_1^2 + q_1 q_2 - q_2 x_1^2 + 
  2 q_1 x_1 x_2 + x_1^2 x_2^2 + q_1 x_1 a_1 - q_2 x_1 a_1\cr 
 &+  q_1 x_2 a_1 + x_1^2 x_2 a_1 + 
  x_1 x_2^2 a_1 + q_1 a_1^2 + x_1 x_2 a_1^2 + q_1 x_1 a_2 \cr 
 &-  q_2 x_1 a_2 + q_1 x_2 a_2 + 
  x_1^2 x_2 a_2 + x_1 x_2^2 a_2 - 
  q_2 a_1 a_2 + x_1^2 a_1 a_2\cr 
 &+  2 x_1 x_2 a_1 a_2 + 
  x_2^2 a_1 a_2 + x_1 a_1^2 a_2 + 
  x_2 a_1^2 a_2 + q_1 a_2^2 +x_1 x_2 a_2^2 \cr 
 &+  x_1 a_1 a_2^2 + x_2 a_1 a_2^2 + a_1^2 a_2^2;}
$$ 
for definition of polynomials $\wt h_m(X_k-Y_l)$ and $\wt e_m(X_k-Y_l)$, 
see Theorems~5 and 6.

\vskip 0.3cm
{\bf (3.4)} Quantum Macdonald polynomials.

\vskip 0.3cm
{\bf Definition 5.} {\it Quantum Macdonald polynomial 
$\wt P_{\ld}(X_{r}~|~p,t)$, $\ld\subset ((n-r)^r)$, is defined from the
decomposition}
$$\wt{\s}_{uw_0}(x,y)=\sum_{\ld\subset ((n-r)^r)}\wt 
P_{\ld}(X_{r}~|~p,t)P_{\wh\ld'}(Y_{n-r}~|~t,p).\eqno (3.9)
$$

\vskip 0.3cm
{\bf Remark.} Polynomial 
$\partial_{w_0^{(n-r)}}^{(y)}\left(\ds\prod_{k=1}^{n-r}
\Delta_{n-k}(y_k~|~X_{n-k})\right)$ is the quantization with respect to 
the $x$ variables of the product 
$\ds\prod_{i=1}^{r}\prod_{j=1}^{n-r}(x_i+y_j)$. Hence, polynomial $\wt 
P_{\ld}(X_r~|~p,t)$ is the quantization of the Macdonald polynomial ~
$P_{\ld}(X_r~|~p,t)$, $\ld\subset ((n-r)^r)$ ([M2], (3.6)).

\vskip 0.5cm
{\bf \S 4. Determinantal formulae for quantum double Schubert polynomials.}
\vskip 0.3cm

In this section we formulate the quantum analogues of determinantal 
formulae for the quantum double and quantum Schubert polynomials 
corresponding to the restricted vexillary 
(see Definition~7), Grassmannian, and 321--avoiding ([BJS], Section~2)
permutations.

\vskip 0.3cm
{\bf (4.1)} Permutations with forbidden subsequences.
\vskip 0.2cm

Let $\tau =(\tau_1,\tau_2\ldots ,\tau_l)\in S_l$ be a permutation.

\vskip 0.3cm
{\bf Definition 6} ([BJS]). {\it A permutation $w\in S_n$ is called 
$\tau$--avoiding if there does not exist $1\le i_1<i_2<\cdots <i_l\le n$ 
such that the following condition holds for all} $1\le j<k\le l$:
$w_{i_j}<w_{i_k}~~~if~ and~ only~ if~~~\tau_j<\tau_k.
$
\vskip 0.2cm

In other words, $w$ has no subsequences in the same relative order as 
$\tau$. We denote by $S_n(\tau )$ the set of all $\tau$--avoiding 
permutations $w\in S_n$. If $\tau_1\in S_{k_1},\ldots ,\tau_r\in 
S_{k_r}$, we denote by $S_n(\tau_1,\tau_2,\ldots ,\tau_r)$ the 
intersection $S_n(\tau_1)\cap\cdots\cap S_n(\tau_r)$. 

Let us consider a few examples.

$\bullet$ A permutation $w$ is 132--avoiding if and only if there do not 
exist $i,j,k$ such that $1\le i<j<k$ and $w_i<w_k<w_j$; it follows, 
from [M1], (1.30), that $w\in S_n(132)$ if and only if $w$ is a dominant 
permutation.

$\bullet$ A permutation $w\in S_n$ is Grassmannian if and only if it 
belongs to the set 
$$S_n(321,2143,3142).
$$ 
The last condition is equivalent 
to the following one: the permutation $w$ has at most one descent.

It is well known ([SS], [BJS]) that $|S_n(132)|=|S_n(321)|=C_n:=\ds{1\over 
n+1}\pmatrix{2n\cr n}$.
\vskip 0.2cm

$\bullet$ A permutation is vexillary if and only if it is 2143--avoiding 
([LS3]; [M1], (1.27)).

\vskip 0.2cm
$\bullet$ A permutation $w\in S_n(2143,1324)$ if and only if the 
Schubert variety $X_w$ is smooth ([LaS]).

\vskip 0.3cm
{\bf Definition 7.} {\it A permutation $w$ is said to be restricted 
vexillary (RV--permutation for short) if it belongs to the set 
$L_n:=S_n(2143,2413,2431)$.}

\vskip 0.3cm
{\bf Remark.} In general, if $w\in L_n$, then inverse permutation 
$w^{-1}$ do not necessary belongs to $L_n$. For example $w=3142\in L_4$, 
but $w^{-1}=2413\not\in L_4$.
It is clear, that a dominant permutation is restricted vexillary, i.e. 
$S_n(132)\subset L_n$. One 
can show that $|S_n(2143)|=|S_n(2431)|=|S_n(2413)|$, and $|L_4|=21$, 
$|L_5|=79$.

\vskip 0.3cm
{\bf (4.2)} Quantum multi--Schur functions.
\vskip 0.2cm

{\bf Definition 8.} {\it Let $X_{k_1},\ldots ,X_{k_n}$ be flagged sets of 
variables and let $\ld ,\mu$ be 
partitions of length $\le n$. The quantum row--flagged Schur function 
$s^q_{\ld /\mu}(X_{k_1},\ldots ,X_{k_n})$ is defined to be
$$s^q_{\ld /\mu}(X_{k_1},\ldots ,X_{k_n})=\det\left(\wt 
h_{\ld_i-\mu_j-i+j}(X_{k_i})\right)_{1\le i,j\le n},\eqno (4.1)
$$
where $\wt h_k(X_r)$ in (4.1) is the quantum complete homogeneous symmetric 
polynomial (see Corollary~2):}
$$\wt h_k(X_r)=\det\left(\wt e_{1-i+j}(X_{r-1+j}\right)_{1\le i,j\le k}.
$$

{\bf Definition 9.} {\it Let ${\cal X}=(X_{k_1},\ldots ,X_{k_m})$ and 
${\cal Y}=(Y_{l_1},\ldots ,Y_{l_m})$ be two families of flagged sets of 
variables and $\ld ,\mu$ be partitions of length $\le m$. The quantum 
multi--Schur function $s_{\ld /\mu}^q({\cal X},{\cal Y})$ is defined to be
$$s_{\ld /\mu}^q({\cal X},{\cal Y})=\det\left(\wt 
h_{\ld_i-\mu_j-i+j}(X_{k_i}-Y_{l_j})\right)_{1\le i,j\le m},\eqno (4.2)
$$
where $\wt h_m(X_k-Y_l):=\ds\sum_{j=0}^m\wt h_{m-j}(X_k)e_j(Y_l)$.}

\vskip 0.2cm
{\bf Definition 10.} {\it Let $X_{l_1},\ldots ,X_{l_m}$ be flagged sets of 
variables and let $\ld ,\mu$ be partitions such that $\ld_1,\mu_1\le m$. 
The quantum column--flagged Schur function $\wt 
s_{\ld'/\mu'}(X_{l_1},\ldots ,X_{l_m})$ is defined to be}
$$\wt s_{\ld'/\mu'}(X_{l_1},\ldots ,X_{l_m})=\det\left(\wt 
e_{\ld'_i-\mu'_j-i+j}(X_{l_j})\right)_{1\le i,j\le m}. \eqno (4.3)
$$

\vskip 0.3cm
{\bf (4.3)} Vexillary permutations ([LS1], (3.1); [M1], (1.27)).
\vskip 0.2cm

Let us remind that a permutation $w$ is vexillary if and only if it is 
2143--avoiding, i.e. there do not exists $i,j,k,l$ such that $1\le 
i<j<k<l$ and $w_j<w_i<w_l<w_k$.

Let $w$ be a permutation with code $c(w)=(c_1,c_2,\ldots )$. For each 
$i\ge 1$ such that $c_i\ne 0$, let
$$g_i=\cases{i, & if for all $j\ge i$, \ $c_j\le c_i$;\cr 
\max\left\{j|j>i\ {\rm and} \ c_j>c_i\right\}, & otherwise .}
$$

Arrange the numbers $g_i$ in increasing order of magnitude, say 
$\theta_1\le\cdots\le\theta_m$. The sequence $\theta 
(w)=(\theta_1,\ldots ,\theta_m)$ is called the flag of $w$. It is a 
sequence of length equal to $l(\ld )$, where $\ld$ is the shape of $w$.

\vskip 0.3cm
{\bf Remark.} The above definition of flag differs from the two 
definitions $\phi (w)$ and $\phi^*(w)$ in [M1], p.14, and the definition 
$\wh\phi (w)$ in [BJS], p.364, or Section~4.4. For example, if 
$w=135624$ is the Grassmannian permutation with code $c=(0,1,2,2)$, then 
$\wh\phi (w)=(234)$, $\theta (w)=(334)$ and $\phi (w)=\phi^*(w)=(444)$. 
The main reason to use the flag $\theta (w)$ instead of flag $\phi :=\phi 
(w)$ is: the formula $\wt{\s}_w=s^q_{\ld}(X_{\phi})$ is valid for smaller 
set of permutations then (4.6). For example, let us take $w=1342$, then 
the code of $w$ is $c=(0,1,1)$ and $\phi =(3,3)$, $\theta =(2,3)$. One 
can check
$$\wt{\s}_{1342}=s^q_{(1,1)}(X_2,X_3)=x_1x_2+x_1x_3+x_2x_3+q_1+q_2,
$$
whereas $s^q_{(1,1)}(X_3,X_3)=x_1x_2+x_1x_3+x_2x_3+q_1+q_2+q_3$.

\vskip 0.3cm
{\bf Proposition 2.} (cf. [LS1], (3.1); [M1], (4.9)). {\it  Let $w$ be 
vexillary 
permutation with shape $\ld =(\ld_1,\ldots ,\ld_m)$ (where $m=l(\ld )$) 
and flag $\phi =(\phi_1,\ldots ,\phi_m)$. Then the Schubert polynomial 
${\s}_w$ is a multi--Schur function, namely}
$${\s}_w=s_{\ld}(X_{\phi_1},\ldots ,X_{\phi_m}).
$$

{\bf Remark.} It follows from Theorem~4, that for the RV--permutations
$$s_{\ld}(X_{\phi})={\s}_w=s_{\ld}(X_{\theta}).
$$

\vskip 0.3cm
{\bf (4.4)} 321--avoiding permutations. 
\vskip 0.2cm

Let us remind that a permutation $w$ is 321--avoiding if it do not 
contains decreasing subsequence of length three. Let $w$ be a 
321--avoiding permutation, follow to [BJS] let us define a skew shape 
$\ld /\mu :=\sigma (w)$ of $w$. Suppose $c(w)=(c_1,\ldots ,c_n)$ and 
$\{\wh\phi_1,\ldots ,\wh\phi_l\} =\{j|c_j>0\}$, with $\wh\phi_1<\cdots 
<\wh\phi_l$. Then $\ld /\mu $ is embedded in ${\bf Z}\times{\bf Z}$ such 
that
$$\ld /\mu =\{(k,h)~|~1\le k\le l, k-\wh\phi_k-c_{\wh\phi_k}+1\le h\le 
k-\wh\phi_k\}.
$$

Let us define the flag $\wh\phi (w)$ of $w$ by
$$\wh\phi (w)=(\wh\phi_1,\ldots ,\wh\phi_l).
$$

{\bf Proposition 3} ([BJS], Theorem 2.2). {\it Let $w$ be a 321--avoiding 
permutation, with skew shape $\sigma (w)=\ld /\mu$ and flag $\wh\phi 
=\wh\phi (w)$. Then
$${\s}_w=s_{\ld /\mu}(X_{\wh\phi}),
$$
the multi--Schur function of shape $\ld /\mu$ and flag $\wh\phi$.}

\vskip 0.3cm
{\bf (4.5)} Determinantal formulae.
\vskip 0.2cm

We start with a generalization of the factorization theorem for Schubert 
polynomials ([M1], (4.6)) to the case of quantum Schubert polynomial 
$\wt{\s}_w$. If $u\in S_m$ and $v\in S_n$, 
let us denote by $u\times v$ the permutation
$u\times v=(u_1,\ldots ,u_m,v_1+m,\ldots ,v_n+m)$ \ in $S_{m+n}$.

\vskip 0.3cm
{\bf Theorem 2.} {\it Let $u\in S_m$ and $v\in S_n$ be permutations. Then}
$$\wt{\s}_{u\times v}=\wt{\s}_u\cdot\wt{\s}_{1^m\times v},\eqno (4.4) 
$$
{\it where $1^m\times v=(1,2,\ldots ,m,v_1+m,\ldots ,v_n+m)$.}

\vskip 0.3cm
{\bf Theorem 3.} {\it Let $w=(w_1,w_2,\ldots ,w_n)\in S_n$ with $w_n=1$. Then
$$\wt{\s}_w=\wt{\s}_u\cdot\wt e_n(X_n),\eqno (4.5)
$$
where $u=(w_1-1,w_2-1,\ldots ,w_{n-1}-1,n)\in S_n$, and $\wt 
e_n(X_n)=\Delta_n(0|X_n)$ is the quantum elementary symmetric function of 
degree $n$ in the variables $X_n=(x_1,\ldots ,x_n)$ (see Section~(2.1)).}

\vskip 0.3cm
{\bf Theorem 4.} {\it Let $w$ be a RV--permutation (see Definition~7) with 
shape $\ld$ and flag $\theta$. Then}
$$\wt{\s}_w=s_{\ld }^q(X_{\theta})=\det\left(\wt 
h_{\ld_i-i+j}(X_{\theta_i})\right)_{1\le i,j\le l(\ld )}. \eqno (4.6)
$$

Theorem~4 follows from more general result:

\vskip 0.3cm
{\bf Theorem 5.} {\it Let $w\in L_n$ be a RV--permutation  with shape 
$\ld =\ld (w)$ and flag $\theta =\theta (w)$ and 
let $\ld'=\ld (w^{-1})$ and $\wt\theta =\theta (w^{-1})$ be shape and 
flag of the inverse permutation $w^{-1}$. Then
$$\eno{
\wt{\s}_w(X_n,Y_n)&=s_{\ld}^q(Z_1,\ldots ,Z_m) &(4.7) \cr
&=\det\left(\wt h_{\ld_i-i+j}(Z_i)\right)_{1\le i,j\le n},}
$$
where $Z_i=X_{\theta_i}-Y_{\wt\theta_{\ld'_i}}$, $1\le i\le m=l(\ld )$, 
and}
~$\wt h_k(X_p-Y_q):=\ds\sum_{j=0}^k\wt h_j(X_p)e_{k-j}(Y_q).$
\vskip 0.3cm

{\bf Remarks.} $i)$ It is well known ([LS1], (3.1); [M1], (1.27)) that 
for a vexillary 
permutation $w$, the shape $\ld (w^{-1})$ of inverse permutation $w^{-1}$ 
coincides with $\ld(w)'$, where $\ld(w)'$ is the conjugate of $\ld (w)$.

$ii)$ It seems that formula (4.7) is new even for ordinary double 
Schubert polynomials.
\vskip 0.3cm

A proof of Theorem~5 is based on the Macdonald method (see [M1], pp.46-49, 
92-93) applied to the variables $Y=(y_1,\ldots ,y_n)$, and the following 
statement which can be checked directly.

\vskip 0.3cm
{\bf Lemma 1.} {\it Let $\delta_n=(n-1,\ldots ,1,0)$ be staircase 
partition. Then
$$s_{\delta_n}^q(Z_1,\ldots ,Z_{n-1})=\Delta_1(y_{n-1}|X_1)
\cdots\Delta_{n-1}(y_1|X_{n-1})=\wt{\s}_{w_0}(x,y),\eqno (4.8)
$$
where $Z_i=X_i-Y_{n-i}$, $1\le i\le n-1$, and $\Delta_k(t|X_k)$ is the 
Givental--Kim determinant (2.2).}

\vskip 0.3cm
{\bf Corollary 6} (of Theorem~6). {\it Let $w\in S_n$ be a dominant 
permutation of shape $\ld$. Then}
$$\wt{\s}_w(x,y)=\det\left(\wt h_{\ld_i-i+j}(X_i-Y_{\ld_i})\right)_{1\le 
i,j\le l(\ld )}.\eqno (4.9)
$$

\vskip 0.3cm
{\bf Theorem 6} (Determinantal formula for quantum factorial Schur functions).

{\it Let $w\in S_n$ be a Grassmannian permutation with 
shape $\ld$ and descent $r$, and let $\wt\theta =\theta (w^{-1})$ be 
the flag of inverse permutation $w^{-1}$. Then
$$\wt{\s}_w(X_n,Y_n)=\det\left(\wt e_{\ld'_i-i+j}(X_{r-1+j}-
Y_{\wt\theta_i})\right)_{1\le i,j\le n-r}, \eqno (4.10)
$$
where  \ $\wt e_m(X_k-Y_l):=\ds\sum_{j=0}^m\wt e_{m-j}(X_k)h_j(Y_l)$.}

\vskip 0.3cm
{\bf Remark.} For general vexillary permutation $w$ with shape $\ld$ and 
flag $\theta$ (or $\phi$, see [M1], p.14) the formula (4.6) (with flag 
$\theta$ or $\phi$) does not valid. For example
$$\eno{
&\wt{\s}_{2431}=\wt s_{\ld}^q(X_{\theta})-q_2q_3, \  \ld =(211), \  
\theta =(223),\cr
&\wt{\s}_{2413}=s_{\ld}^q(X_{\theta})+q_2(x_1+x_2+x_3), \  \ld =(21), \ 
\theta =(2,2).}
$$

It seems plausible, that if $w\in S_n(2143)$, but $w\not\in L_n$, then 
$\wt{\s}_w\ne s_{\ld}^q(X_{\theta})$. For example, let us take the 
vexillary permutation $w=42513\not\in L_5$. Then $c=(3120)$, $\ld 
=(321)$, $\theta =(133)$, and
$$s^q_{\ld}(X_{\theta})=\wt{\s}_w+q_3\wt{\s}_{41235}\cdot\wt{\s}_{12354}
-q_3\wt{\s}_{51234}.
$$
Hence, \ \ $s_{\ld}^q(X_{\theta})\vert_{q_3=0}=\wt{\s}_w\vert_{q_3=0}$. It 
is easy to see that
$${\partial\over\partial 
q_3}\wt{\s}_{42513}=-\wt{\s}_{42135}.
$$

{\bf Exercise.} Let $w\in S_n$ be a vexillary permutation with code 
$c=(c_1,\ldots ,c_{n-1})$, shape $\ld$, and flag $\theta$. Let us define 
$m=\max\{j|c_j\ne 0\}$. Prove
$$s_{\ld}^q(X_{\theta})|_{\matrix{{\scriptstyle q_m=0,}\cr
{\scriptstyle q_{m+1}=0,}\cr\cdots}}=
\wt{\s}_w(x)|_{\matrix{{\scriptstyle q_m=0,}\cr{\scriptstyle q_{m+1} =0,}
\cr \cdots}}.
$$

Note, that if \hbox{$w\in L_n$} be a random permutation with code 
$c=(c_1,\ldots ,c_{n-1})$ and \hbox{$m=\max\{ j|c_j\ne 0\}$,} 
then \ $\ds{\partial\over
\partial q_m}\wt{\s}_w\ne 0$.

\vskip 0.3cm
{\bf Problem.} To describe the matrix elements of the operators 
$\ds{\partial\over\partial q_i}$, and $(x_1+\cdots 
+x_i)\ds{\partial\over\partial q_i}$ in the basis of quantum Schubert 
polynomials.
\vskip 0.3cm

{\bf (4.6)} Quantum double Schubert polynomials for 321--avoiding 
permutations.
\vskip 0.3cm

{\bf Conjecture 1.} {\it Let $w$ be a 321--avoiding permutation of shape 
$\sigma (w)=\ld /\mu$ and flag $\wh\phi :=\wh\phi (w)$. Then}
$$\wt{\s}_w=s_{\ld /\mu}^q(X_{\wh\phi}).
$$

{\bf Conjecture 2.} {\it 
More generally, let $w\in S_n$ be a 321--avoiding permutation with shape 
$\sigma (w)=\ld /\mu$ and flag $\wh\phi =\wh\phi (w)$ and let 
$\wh\phi'=\wh\phi (w^{-1})$ be the flag of  inverse permutation $w^{-1}$. 
Then
$$\wt{\s}_w(X_n,Y_n)=\det\left(\wt h_{\ld_i-\mu_j-i+j}(X_{\wh\phi_i}-
Y_{\wh\phi'(j)})\right)_{1\le i,j \le l(\ld )}, \eqno (4.11)
$$
where $\wh\phi'(j)=\wh\phi'_{\sigma (w^{-1})_j}$, and $\sigma (w^{-1})_j$ 
is the length of $j$-th row in (skew) diagram $\sigma (w^{-1})$.}
\vskip 0.3cm

Conjecture~2 suggests the interesting identities between quantum 
row--flagged and quantum column--flagged Schur functions. For example, 
let us take the Grassmannian permutation $w=2413$ (note that $w\not\in 
L_4$). Then $c(w)=(1,2)$, $\ld (w)=(2,1)$, $\theta (w)=(2,2)$, 
$c(w^{-1})=(2,0,1)$, $\ld (w^{-1})=(2,1)$, $\theta (w^{-1})=(1,3)$. 
According to the Theorem~6,
$$\wt{\s}_{2413}(x,y)=\left|\matrix{\wt e_2(X_2-Y_1)&\wt e_3(X_3-Y_1)\cr 
1&\wt e_1(X_3-Y_3)}\right|.
$$
But permutation 2413 is also 321--avoiding with shape $\sigma 
(w)=(2,2)/(1,0)$, flag $\wh\phi =(1,2)$, and inverse flag 
$\wh\phi'=(1,3)$. One can check that (cf. Conjecture~2)
$$\wt{\s}_{2413}(x,y)=\left|\matrix{\wt h_1(X_1-Y_1)&\wt h_3(X_1-Y_3)\cr 
1&\wt h_2(X_2-Y_3)}\right|.
$$

\vskip 0.3cm
{\bf (4.7)} Grassmannian permutations and duality theorem.
\vskip 0.2cm

Let us consider the case of Grassmannian permutations in more details. 
Note at first that each Grassmannian permutation is also 321--avoiding.
Let $w\in S_n$ be Grassmannian permutation with shape $\ld 
=(m^{k_m},(m-1)^{k_{m-1}},\ldots ,2^{k_2},1^{k_1})$ and descent $r$ (we 
assume that $n\ge r+m$). Then the code and flag of $w$ are
$$\eno{
&c(w)=(\underbrace{0,\ldots ,0}_{k_0},\underbrace{1,\ldots ,1}_{k_1},
\ldots ,\underbrace{m,\ldots ,m}_{k_m},0,\ldots ,0),\cr
&\theta (w)=(k_0+1,k_0+2,k_0+3,\ldots ,k_0+k_1+\cdots +k_m).}
$$
It is clear that $l:=l(\ld )=k_1+\cdots +k_m$, $r=k_0+l(\ld )$, and the 
shape of $w$ as a 321--avoiding permutation is a skew shape $\sigma 
(w)=(\ld_1^l)/(\ld_1-\ld_l,\ld_1-\ld_{l-1},\ldots ,\ld_1-\ld_2)$ (see 
[BJS], (11), or Section (4.4)), and $\wh\phi (w)=\theta (w)$. Now let us 
consider the inverse permutation $w^{-1}\in L_n$. Then (cf. [M1], (1.25))
$$c(w^{-1})=(\underbrace{0,\ldots ,0}_{k_0},\ld_1',\underbrace{0,\ldots 
,0}_{k_1},\ld_2',\underbrace{0,\ldots ,0}_{k_2}\ld_3',\ldots 
,\ld_{m-1}',\underbrace{0,\ldots ,0}_{k_{m-1}},\ld_m', 0,\ldots ,0),
$$
where $\ld'=(\ld_1',\ldots ,\ld_m')$ is the conjugate of $\ld$. Finally, 
the flag of inverse permutation is $\wt\theta =(\theta 
(w^{-1})_i)$, where $\theta (w^{-1})_i=\ld_1'-\ld_i'+i$, $1\le i\le m$. 

The following result is a partial confirmation of Conjecture~2:
\vskip 0.3cm

{\bf Theorem 7.} {\it Let \ $\ld =(m^{k_m},(m-1)^{k_{m-1}},\ldots 
,2^{k_2},1^{k_1})$ \ be a partition of length\break $r=k_1+\cdots +k_m$. Then
$$\det\left(\wt e_{\ld'_i-i+j}(X_{r-1+j}-Y_{r-\ld_i'+i})\right)_{1\le i,j\le m}
=\det\left(\wt h_{\ld_i-i+j}(X_{r-j+1}-Y_{\gamma (r-i+1)})\right)_{1\le i,j\le 
r},
$$
where $\gamma (j)=r+\ld_j-\ld_{\ld_j}'$.}

\vskip 0.2cm
Proof of Theorem~7 follows from the following relations between quantum 
elementary symmetric functions $\wt e_k(X_{n+m-1})$ and quantum complete 
homogeneous functions $\wt h_j(X_n)$:

\vskip 0.3cm
{\bf Lemma 2.} ~{\it If $n\ge 1$, then}
$$\ds\sum_{j=0}^m(-1)^j\wt e_{m-j}(X_{n+m-1})\wt 
h_j(X_n)=\delta_{m,0}.
$$

\vskip 0.3cm
{\bf Corollary 7} (Jacobi--Trudi formula for quantum Schur functions). 

{\it Let $\ld$ be a partition, $l(\ld )\le n$. Then}
$$\wt s_{\ld}(X_n)=\det\left(\wt h_{\ld_i-i+j}(X_{n+1-j})
\right)_{1\le i,j\le n}.
$$
\vskip 0.2cm 

Theorem~7 is a special case of general duality theorem for quantum 
multi--Schur functions (duality theorem for ordinary multi--Schur 
functions see in [M1], (3.8'), and [M2], (8.3)). We suppose to 
consider the general case in a separate publication.

\vskip 0.3cm
{\bf (4.8)} Stable quantum Schubert polynomials.
\vskip 0.2cm

Let $w\in S_n$ be a permutation, consider the permutations 
$$w_m:=1^m\times w=(1,\ldots ,m,w_1+m,\ldots ,w_n+m)\in S_{n+m},~~ 
m=0,1,2,\ldots.
$$ 
One can show that there exists the limit
$$\wt F_w:=\lim_{m\to\infty}\wt{\s}_w\in{\bf Z}[x_1,x_2,\ldots 
;q_1,q_2,\ldots ].
$$
The function $\wt F_w$ is called to be the quantum stable Schubert 
polynomial corresponding to the permutation $w$. In classical case 
$q_1=q_2=\cdots =0$, the functions~ $F_w:=\wt F_w|_{q_1=q_2=\cdots =0}$, 
were studied in [LS3]; [S]; [M1], Chapter~VII; ...

\vskip 0.3cm
{\bf Example.} Let us take $w=321\in S_3$. Then the shape and flag of 
permutation $w_m=1^m\times w$ are $\ld =(21)$ and $\theta =(m+1,m+2)$. 
Hence,
$$\wt{\s}_{w_m}=\left|\matrix{\wt h_2(X_{m+1})&\wt h_3(X_{m+1})\cr
1&\wt h_1(X_{m+2})}\right|,
$$
and
$$\wt F_{21}:=\lim_{m\to\infty}\wt{\s}_{w_m}=\left|\matrix{\wt 
h_2(X_{\infty})& \wt h_3(X_{\infty})\cr 1 &\wt h_1(X_{\infty})}\right|
=s_{21}(x)+\sum_{i\ge 1}q_i(x_i+x_{i+1})=\wt s_{21}(X_{\infty}).
$$
Note, that the function $\wt F_{21}$ does not symmetric with respect to the 
$x$ variables.

It seems an interesting problem to understand the meaning of quantum 
quasi--symmet\-ric functions $\wt Q_{D,p}(X_m)$, which are the quantization 
of the fundamental quasi--symmetric functions $Q_{D,p}(x_1,\ldots ,x_m)$ 
([S], p.360; [M1], p.101).

\vskip 0.3cm
{\bf (4.9)} What is next?
\vskip 0.2cm

$1^0$  Let us call a permutation $w$, for which the quantum Schubert 
polynomial $\wt{\s}_w$ is a quantum row--flagged skew Schur function 
$s^q_{\ld /\mu}(X_{\theta})$, a restricted skew vexillary permutation 
(RSV--permutation for short) of shape $\ld /\mu$. It is an interesting 
problem to classify the RSV--permutations.

$2^0$ It seems very interesting to find a quantum analog of the Magyar 
algebra--geometric approach ([Ma1]) to the theory of Schubert polynomials.

$3^0$ What is the group representation meaning of the quantum Schur 
functions?

$4^0$ What is a quantum analog of the Macdonald operators ([M3])? In other 
words, the eigenfunctions of which commuting family of operators the 
quantum Macdonald polynomials are?

$5^0$  Part of our results can be generalized to the case of quantum 
double Grothendieck polynomials [K], and quantum Key polynomials. The 
work is in progress and we hope to present our results in the nearest 
future. 

\vfil\eject
{\bf References.}
\vskip 0.5cm

\item{[BJS]} Billey S., Jockusch W. and Stanley R., {\it Some 
combinatorial properties of Schubert polynomials,} Journ. Algebraic Comb.,
1993, v.2, p.345-374;

\item{[C]} Ciocan--Fontanine I., {\it Quantum cohomology of flag varieties},
        Intern. Math. Research Notes, 1995, n.6, p.263-277;

\item{[FGP]}  Fomin S, Gelfand S. and Postnikov A., {\it Quantum Schubert
     polynomials,} Preprint\break\hbox{AMSPPS} \#199605--14--008, 1996, 44p.;

\item{[GK]} Givental A., Kim B., {\it Quantum cohomology of flag manifolds and 
        Toda lattices}, Comm. Math. Phys., 1995, v.168, p.609-641;

\item{[K]} Kirillov A.N., {\it Quantum Grothendieck polynomials,} 
          Preprint, 1996, q-alg/9610034, 12p;        
        
\item{[KM]} Kirillov A.N. and Maeno T., {\it Quantum double Schubert 
 polynomials, quantum Schubert polynomials and Vafa--Intriligator formula,}
 Preprint, 1996, q-alg/9610022, 52p.;
 
\item{[L]} Lascoux A., {\it Schubert polynomials and shifted Schur 
functions,} Manuscript, 1995.

\item{[LS1]} Lascoux A. and Sch\"utzenberger M.-P., {\it Polynomes de Schubert,}
C.R. Acad. Sci. Paris, 1982, v.294, p.447-450;

\item{[LS2]} Lascoux A. and Sch\"utzenberger M.-P., {\it Symmetry and flag 
manifolds,} Lect. Notes in Math., 1983, v.996, p.118-144;

\item{[LS3]} Lascoux A. and Sch\"utzenberger M.-P., {\it Structure de 
Hopf de l'anneu de cohomologie et de l'anneau de Grothendieck d'une 
variete de drapeaux,} C.R. Acad. Sci. Paris, 1982. v.295, p.629-633;

\item{[LaS]} Lakshmibai V. and Sandhya B., {\it Criterion for smoothness 
of Schubert varieties,} Proc. Indian Acad. Sci. M., 1990, v.100, p.45-92;       

\item{[M1]} Macdonald I.G., {\it Notes on Schubert polynomials}, Publ.
 LCIM, 1991, Univ. of Quebec a Montreal;

\item{[M2]} Macdonald I.G., {\it Schur functions: theme and variations,} 
Publ. I.R.M.A. Strasbourg, 1992, 498/S-28, Actes 28-e Seminaire Lotharingien, 
p.5--39;

\item{[M3]} Macdonald I.G., {\it Symmetric functions and Hall 
polynomials,} Second ed., Oxford Univ. Press, New York/London, 1995;

\item{[Ma]} Magyar P., {\it Bott--Samelson varieties and configuration 
spaces,} Preprint, 1996, alg-geom/9611019, 32p.;

\item{[S]} Stanley R.P., {\it On the number of reduced decompositions of 
elements of Coxeter groups,} European J. Comb., 1984, v.5, p359-372;

\item{[SS]} Simion R. and Schmidt F., {\it Restricted permutations,} 
European J. Comb., 1985, v.6, p.383-406.

\end